\documentclass[11pt]{amsart}
\usepackage{amscd, amsfonts, amsthm, amssymb}
\title{The groups of PL and Lipschitz homeomorphisms of noncompact 2-manifolds} 
\author{Tatsuhiko Yagasaki} 
\address{Department of Mathematics,  Kyoto Institute of Technology,  Matsugasaki,  Sakyoku, Kyoto 606, Japan}
\email{yagasaki@ipc.kit.ac.jp}

\newtheorem{theorem}{Theorem}[section]
\newtheorem{proposition}{Proposition}[section] 

\newtheorem{lemma}{Lemma}[section]
\newtheorem*{claim}{Claim}

\theoremstyle{definition}

\newtheorem{problem}{Problem}
\newtheorem{example}{Example}[section]

\newtheorem*{Proof of Proposition 3.1}{Proof of Proposition 3.1}

\def \cal {\mathcal}
\def \phi {\varphi}

\begin{document}

\baselineskip    6 mm
\maketitle
\begin{abstract} 
Suppose $M$ is a noncompact connected PL 2-manifold. 
In this paper we study the topological property of the triple 
$({\cal H}(M)_0$, ${\cal H}^{{\rm PL}}(M)_0, {\cal H}^{{\rm PL, c}}(M)_0)$, 
where ${\cal H}(M)_0$ is the identity component of the homeomorphism group ${\cal H}(M)$ of $M$ with the compact-open topology,
and ${\cal H}^{{\rm PL}}(M)_0$ and ${\cal H}^{{\rm PL, c}}(M)_0$ are 
the identity components of the subgroups consisting of PL-homeomorphisms of $M$ and ones with compact supports. 
We show that this triple is a $(s^{\infty},\sigma^{\infty},\sigma^{\infty}_f)$-manifold and determine its topological type. 
We also study the subgroups of Lipschitz homeomorphisms. 
\end{abstract}

\section{Introduction}
This article is a continuation of study of homeomorphism groups of 2-manifolds and related topics \cite{Ya1, Ya2, Ya3, Ya4, Ya5}. 
The purpose of this paper is to study topological properties of the subgroups of PL and Lipschitz homeomorphisms of 
2-manifolds with the compact-open topology \cite{Ya2}. 
The property of groups of quasiconformal homeomorphisms on Riemann surfaces are studied in \cite{Ya5}. 

For a manifold $M$ and a subset $X$ of $M$, let ${\cal H}_X(M)$ denote 
the group of homeomorphisms $h$ of $M$ onto itself with $h|_X =id$, equipped with the compact-open topology. 
When $M$ is a PL-manifold, ${\cal H}_X^{\rm PL}(M)$ denotes 
the subgroup of ${\cal H}_X(M)$ consisting of PL-homeomorphisms, 
and when $M$ has a fixed metric, ${\cal H}_X^{\rm (loc) \, LIP}(M)$ denotes 
the subgroup of (locally) Lipschitz homeomorphisms. 
In addition, the superscript \lq\lq$c$" means the compact supports and the subscript \lq\lq 0" means 
the identity connected component.  

When $M$ is a compact PL 2-manifold and $X$ is a subpolyhedron of $M$, 
the pair $({\cal H}_X(M), {\cal H}_X^{{\rm PL}}(M))$ is a $(s, \sigma)$-manifold \cite{GH}. 
(See Section 2 for the basic notions on infinite dimensional topological manifolds.) 
We extend this result to the noncompact case. 
Suppose $M$ is a noncompact connected PL 2-manifold and $X$ is a compact subpolyhedron of $M$. 
We have already shown that ${\cal H}_X(M)_0$ is a $s$-manifold and determined its topological type \cite{Ya4}:

\renewcommand{\labelenumi}{(\Roman{enumi})}
\begin{enumerate}
\item ${\cal H}_X(M)_0 \cong s \times {\Bbb S}^1$ if $(M, X) \cong ({\Bbb R}^2, \emptyset)$, $({\Bbb R}^2, 1pt)$, 
$({\Bbb S}^1 \times {\Bbb R}^1, \emptyset)$, $({\Bbb S}^1 \times [0, 1), \emptyset)$ or $({\Bbb P}^2 \setminus 1pt, \emptyset)$
\item ${\cal H}_X(M)_0 \cong s$ in all other cases, 
\end{enumerate}
where ${\Bbb R}^n$ is the Euclidean $n$-space, ${\Bbb S}^n$ is the $n$-sphere and ${\Bbb P}^2$ is the projective plane.

In this article we show that ${\cal H}_X^{{\rm PL, c}}(M)_0$ has the homotopy negligible complement in ${\cal H}_X(M)_0$. 
By the characterization and homotopy invariance of infinite dimensional manifolds we obtain the following conclusion:

\begin{theorem}
If $M$ is a noncompact connected PL 2-manifold and $X$ is a compact subpolyhedron of $M$, 
then $({\cal H}_X(M)_0, {\cal H}_X^{{\rm PL}}(M)_0, {\cal H}_X^{{\rm PL,c}}(M)_0)$ is 
a $(s^{\infty},\sigma^{\infty},\sigma^{\infty}_f)$-manifold. 
This triple is homeomorphic to $(s^{\infty},\sigma^{\infty},\sigma^{\infty}_f) \times {\Bbb S}^1 $ in the case (I) 
and to $(s^{\infty},\sigma^{\infty},\sigma^{\infty}_f)$ in the case (II). 
\end{theorem}

Next we incorporate subgroups of Lipschitz homeomorphisms. 
A Euclidean PL-manifold means a (closed) subpolyhedron of a Euclidean space ${\Bbb R}^N$ which is a PL-manifold with respect to this triangulation 
and is equipped with the metric induced from the standard metric of ${\Bbb R}^N$.
We have shown that when $M$ is a compact Euclidean PL 2-manifold and $X$ is a subpolyhedron of $M$, 
the triple $({\cal H}_X(M), {\cal H}_X^{\rm LIP}(M), {\cal H}_X^{\rm PL}(M))$ is a $(s, \Sigma, \sigma)$-manifold \cite{Ya1}. 
The following is a noncompact version: 

\begin{theorem}
If $M$ is a noncompact connected Euclidean PL 2-manifold and $X$ is a compact subpolyhedron of $M$, then the quadruples 
$({\cal H}_X(M)_0, {\cal H}_X^{\rm loc \, LIP}(M)_0, {\cal H}_X^{\rm LIP, c}(M)_0, {\cal H}_X^{\rm PL, c}(M)_0)$ and \\
$({\cal H}_X(M)_0, {\cal H}_X^{\rm loc \, LIP}(M)_0, {\cal H}_X^{\rm LIP}(M)_0, {\cal H}_X^{\rm PL, c}(M)_0)$ are 
$(s^{\infty},\Sigma^{\infty},\Sigma^{\infty}_f, \sigma^{\infty}_f)$-manifolds. 
These quadruples are homeomorphic to $(s^{\infty},\Sigma^{\infty},\Sigma^{\infty}_f, \sigma^{\infty}_f) \times {\Bbb S}^1$ 
in the case (I) and to $(s^{\infty},\Sigma^{\infty},\Sigma^{\infty}_f, \sigma^{\infty}_f)$ in the case (II). 
\end{theorem}

Any smooth $n$-manifold $M$ has a $C^1$-triangulation \cite{Wh} and it can be embedded into 
the $(2n+1)$-dimensional Euclidean space as a subpolyhedron. 
Since this embedding is locally Lipschitz for any Riemannian metric on $M$, 
Theorem 1.2 implies the corresponding statements on Riemannian 2-manifolds. 
Since this embedding is not necessarily Lipschitz, 
to obtain a statement on subgroups of Lipschitz homeomorphisms on noncompact Riemannian manifolds, 
we need to modify the proof itself of Theorem 1.2. 

\begin{theorem} 
Suppose $M$ is a connected Riemannian 2-manifold with a $C^1$-triangulation and $X$ is a compact subpolyhedron of $M$. \\
(i) If $M$ is compact and $X \neq M$, then $({\cal H}_X(M), {\cal H}_X^{\rm LIP}(M), {\cal H}_X^{\rm PL}(M))$ is 
a $(s, \Sigma, \sigma)$-manifold. \\
(ii) If $M$ is noncompact, then 
$({\cal H}_X(M)_0, {\cal H}_X^{\rm loc \, LIP}(M)_0, {\cal H}_X^{\rm LIP, c}(M)_0, {\cal H}_X^{\rm PL, c}(M)_0)$ 
is a $(s^{\infty},\Sigma^{\infty},\Sigma^{\infty}_f, \sigma^{\infty}_f)$-manifold. \\
(iii) If $M$ is a noncompact complete Riemannian manifold, then 
$({\cal H}_X(M)_0, {\cal H}_X^{\rm LIP}(M)_0, {\cal H}_X^{\rm PL, c}(M)_0)$ is a $(s, \Sigma, \sigma)$-manifold and 
$({\cal H}_X(M)_0, {\cal H}_X^{\rm loc \, LIP}(M)_0, {\cal H}_X^{\rm LIP}(M)_0)$ is 
a $(s^{\infty},\Sigma^{\infty},\Sigma^{\infty}_f)$-manifold. 
\end{theorem}


This paper is organized as follows: 
Section 2 contains some basic facts on infinite-dimensional topological manifolds and a characterization of 
$(s^{\infty},\Sigma^{\infty},\Sigma^{\infty}_f, \sigma^{\infty}_f)$-manifolds. 
Its proof is included in Section 5. 
In Sections 3 and 4 we show Homotopy absorption property, Class property and 
Stability property of homeomorphism groups.
Combined with the characterization of $(s^{\infty},\Sigma^{\infty},\Sigma^{\infty}_f, \sigma^{\infty}_f)$-manifolds, 
these properties imply Theorems 1.1 -- 1.3. 
Throughout the paper spaces are assumed to be separable, metrizable and maps to be continuous. 
By ${\cal C}(X, Y)$ we denote the space of maps from a space $X$ to a space $Y$ with the compact-open topology. 
For a subset $A$ of a metric space $X$, $N_X(A, \varepsilon)$ denotes the $\varepsilon$-neighborhood of $A$ in $X$.

\section{Characterization of infinite-dimensional topological manifolds}
\subsection{Basic facts on infinite-dimensional manifolds} \mbox{} 

In this subsection we collect basic facts on infinite-dimensional manifolds. 
As for the model spaces of infinite-dimensional manifolds we follow the standard convension: 
$s = {\Bbb R}^{\infty} \ (\cong \ell_2)$, 
$\sigma = \{(x_n) \in s \, : \, x_n = 0$ (almost all $n)\}$, 
$\Sigma = \{(x_n) \in s \, : \, \sup_n\,|x_n| < \infty \}$, 
$\sigma^{\infty}_f = \{(x_n) \in \sigma^{\infty} \, : \, x_n = 0$ (almost all $n)\}$, 
$\Sigma^{\infty}_f = \{(x_n) \in \Sigma^{\infty} \, : \, x_n = 0$ (almost all $n)\}$ and 
$s^{\infty}_b = \{ (x_n) \in s^{\infty} : \, \sup_{n,i} |(x_n)_i| < \infty \}$. 
We note that (a) $(s^{\infty}, \Sigma^{\infty}_f, \sigma^{\infty}_f) \cong (s, \Sigma, \sigma)$ \cite[Theorem 1.1]{Ya1}, 
(b) $(s^{\infty}, \Sigma^{\infty}, \sigma^{\infty}_f) \cong (s^{\infty}, \sigma^{\infty}, \sigma^{\infty}_f)$ \cite[Corollary 3.16]{Ya1}, and 
(c) $(s^{\infty}, \Sigma^{\infty}, \Sigma^{\infty}_f, \sigma^{\infty}_f) 
\cong (s^{\infty}, \Sigma^{\infty}, s^{\infty}_b, \sigma^{\infty}_f)$ (\S5.3. Proposition 5.2). 

We say that a subset $Y$ has the homotopy negligible (h.n.) complement in $X$ 
if there exists a homotopy $\phi_t : X \to X \ (0 \leq t \leq 1)$ such that $\phi_0 = id_X$ 
and $\phi_t(X) \subset Y \ (0 < t \leq 1)$. 
The homotopy $\phi_t$ is called an absorbing homotopy of $X$ into $Y$. 
A space is $\sigma$-(fd-) compact if it is a countable union of (finite dimensional) compact subsets. 
A tuple $(X, X_1, \cdots, X_n)$ means a tuple of a space $X$ and its subspaces $X_1 \supset \cdots \supset X_n$. 
Suppose $(E, E_1, \cdots, E_n)$ is a model tuple. A tuple $(X, X_1, \cdots, X_n)$ is said to be $(E, E_1, \cdots, E_n)$-stable 
if $(X \times E, X_1 \times E_1, \cdots, X_n \times E_n) \cong (X, X_1, \cdots, X_n)$ (a homeomorphism of tuples). 
We say that a tuple $(X, X_1, \cdots, X_n)$ is a $(E, E_1, \cdots, E_n)$-manifold 
if each point $x$ in $X$ has a neighborhood $U$ such that $(U, U \cap X_1, \cdots, U \cap X_n) \cong (E, E_1, \cdots, E_n)$. 
In this article we apply the following characterization of 
$(s^{\infty}, \Sigma^{\infty}, \Sigma^{\infty}_f, \sigma^{\infty}_f)$-manifolds: 

\renewcommand{\labelenumi}{(\roman{enumi})}
\begin{theorem}
A quadruple $(X, X_1, X_2, X_3)$ is a $(s^{\infty}, \Sigma^{\infty}, \Sigma^{\infty}_f, \sigma^{\infty}_f)$-manifold iff
\vspace{-5pt}
\begin{enumerate}
\setlength{\itemsep}{-2pt}
\item $X$ is a separable completely metrizable ANR, 
\item (Homotopy absorption property) $X_3$ has the h.n.\ complement in $X$, 
\item (Class condition) $X_1$ is $F_{\sigma \delta}$ in $X$, $X_2$ is $\sigma$-compact and $X_3$ is $\sigma$-fd-compact, 
\item (Stability) $(X, X_1, X_2, X_3)$ is $(s^{\infty}, \Sigma^{\infty}, \Sigma^{\infty}_f, \sigma^{\infty}_f)$-stable. 
\end{enumerate}
\end{theorem}

The topological types of these manifolds are detected by their homotopy types. 

\begin{proposition} {\rm (Homotopy invariance)}  
Suppose $(X, X_1, X_2, X_3)$ and $(Y, Y_1, Y_2, Y_3)$ are \\
$(s^{\infty}, \Sigma^{\infty}, \Sigma^{\infty}_f, \sigma^{\infty}_f)$-manifolds. 
Then $(X, X_1, X_2, X_3) \cong (Y, Y_1, Y_2. Y_3)$ iff $X \simeq Y$ {\rm (}homotopy equivalent{\rm )}. 
In particular, if $X$ has a homotopy type of a compact polyhedron $P$, then $(X, X_1, X_2, X_3) 
\cong (s^{\infty}, \Sigma^{\infty}, \Sigma^{\infty}_f, \sigma^{\infty}_f) \times P$.
\end{proposition}

In the above assertions, if we omit the conditions related to $\Sigma^{\infty}$ (respectively $\Sigma^{\infty}_f$ or $\sigma^{\infty}_f$) 
then we obtain the characterization and homotopy invariance of manifolds modeled on 
$(s, \Sigma, \sigma) \cong (s^{\infty}, \Sigma^{\infty}_f, \sigma^{\infty}_f)$ 
(respectively $(s^{\infty}, \sigma^{\infty}, \sigma^{\infty}_f) \cong (s^{\infty}, \Sigma^{\infty}, \sigma^{\infty}_f)$ 
or $(s^{\infty}, \Sigma^{\infty}, \Sigma^{\infty}_f)$).
Theorem 2.1 and Proposition 2.1 will be proved in Section 5.  

\subsection{Proof of Theorems 1.1 -- 1.3} \mbox{} 

To prove Theorems 1.1 -- 1.3 we have to check the conditions (i) -- (iv) in Theorem 2.1.
As for the condition (i), if $X$ is a locally compact, locally connected, separable metrizable space and $A$ is a closed subset of $X$, 
then ${\mathcal H}_A(X)$ is a separable completely metrizable topological group. 
The ANR property follows from \cite{LM}, \cite[Corollary 1.1]{Ya4}: 

\begin{proposition} Suppose $M$ is a PL $2$-manifold and $X$ is a compact subpolyhedron of $M$. 
If $M$ is compact then ${\cal H}_X(M)$ is an ANR, and if $M$ is noncompact and connected then ${\cal H}_X(M)_0$ is an ANR.    
\end{proposition}

The conditions (ii), (iii) and (iv) will be verified in Sections 3, 4.1 and 4.2 respectively. 
As mentioned in Section 1, in the noncompact case the homotopy type of ${\cal H}_X(M)_0$ is classified as follows \cite[Theorem 1.1]{Ya4}:

\begin{proposition} 
Suppose $M$ is a noncompact connected PL 2-manifold and $X$ is a compact subpolyhedron of $M$. 
\renewcommand{\labelenumi}{(\Roman{enumi})}
\begin{enumerate}
\item ${\cal H}_X(M)_0 \simeq {\Bbb S}^1$ if $(M, X) \cong 
({\Bbb R}^2, \emptyset)$, $({\Bbb R}^2, 1pt)$, $({\Bbb S}^1 \times {\Bbb R}^1, \emptyset)$, 
$({\Bbb S}^1 \times [0, 1), \emptyset)$ or $({\Bbb P}^2 \setminus 1pt, \emptyset)$,
\item ${\cal H}_X(M)_0 \simeq \ast$ in all other cases.
\end{enumerate}
\end{proposition}

\noindent The assertions on topological type in Theorems 1.1 and 1.2 follow from Propositions 2.1 and 2.3.

In Section 3 we use the fact that the restriction map from a homeomoprhism group to an embeddding space is a principal bundle. 
We conclude this section with a statement on this fact \cite[Corollary 1.1]{Ya3}: 
Suppose $M$ is a PL 2-manifold and 
$K \subset X$ are compact subpolyhedra of $M$. Let ${\cal E}_K(X, M)$ denote 
the space of embeddings $f : X \hookrightarrow M$ with $f|_K = id$, equipped with the compact-open topology. 
We consider the subspace ${\cal E}_K(X, M)^{\ast} = \{f \in {\cal E}_K(X,M) \, 
: \, f(X \cap \partial M) \subset \partial M, f(X \cap {\rm Int} \,M) \subset {\rm Int} \, M \}$ 
(the space of proper embeddings). Let ${\cal E}_K(X, M)^{\ast}_0$ denote the connected component of 
the inclusion $i_X : X \subset M$ in ${\cal E}_K(X, M)^{\ast}$. 

\begin{proposition}
For any open neighborhood $U$ of $X$ in $M$, the restriction map 
$\pi : {\cal H}_{K \cup (M \setminus U)}(M)_0 \to {\cal E}_K(X, U)_0^{\ast}$, $\pi(f) = f|_X$, 
is a principal bundle with fiber ${\cal G} \equiv {\cal H}_{K \cup (M \setminus U)}(M)_0 \cap {\mathcal H}_X(M)$, 
where the group ${\mathcal G}$ acts on ${\cal H}_{K \cup (M \setminus U)}(M)_0$ by right composition. 
\end{proposition}

\section{Homotopy absorption property of ${\cal H}_X^{{\rm PL}}(M)$}

The purpose of this section is to prove the next assertion:

\begin{proposition} 
Suppose $M$ is a PL 2-manifold and $X$ is a compact subpolyhedron of $M$. If $M$ is compact then ${\cal H}_X^{{\rm PL}}(M)$ 
has the h.n.\ complement in ${\cal H}_X(M)$, and 
if $M$ is noncompact connected then ${\cal H}_X^{{\rm PL, c}}(M)_0$ has the h.n.\ complement in ${\cal H}_X(M)_0$. 
\end{proposition} 

We proceed to the verification of Proposition 3.1. 
The compact case is already known \cite{GH},cf.\,\cite[Fact 4.3]{Ya3}. Hence below we assume that $M$ is noncompact and connected. 
We use the same notations as in \cite[Section 4]{Ya3}: We set $M_0 = X$ and write as $M = \cup_{i=0}^{\infty} \, M_i$, 
where for each $i \geq 1$ (a) $M_i$ is a nonempty compact connected PL 2-submanifold of $M$ and $M_{i-1} \subset {\rm Int}_M M_i$, 
(b) for each component $L$ of $cl \, (M \setminus M_i)$, $L$ is noncompact and $L \cap M_{i+1}$ is connected, 
and (c) $M_1 \cap \partial M \neq \emptyset$ if $\partial M \neq \emptyset$. 
Taking a subsequence, we have the following cases: (i) each $M_i$ is a disk, 
(ii) each $M_i$ is an annulus (the inclusion $M_i \subset M_{i+1}$ is essential), 
(iii) each $M_i$ is a M\"obius band, and (iv) each $M_i$ is not a disk, an annulus or a M\"obius band. 
We choose a metric $d$ on $M$ with $d \leq 1$ and metrize ${\cal H}_X(M)$ by the metric $\rho$ defined by 
\[ \rho(f, g) = \sum_{i=1}^{\infty} \, \frac{1}{2^i} \sup_{x \in M_i} \, d(f(x), g(x)). \]

We also use the following notations: 
Let $U_j = {\rm Int}_M M_j$ and $L_j = {\rm Fr}_M M_j$ ($j \geq 1$). 
For $j > i \geq 0$ let ${\cal H}_j = {\cal H}_{X \cup (M \setminus U_j)}(M)_0$, 
${\cal U}^i_j = {\cal E}_X(M_i, U_j)^{\ast}_0$, and let $\pi^i_j : {\cal H}_j \to {\cal U}^i_j$ and $\pi^i : {\cal H}_X(M)_0 \to {\cal E}_X(M_i, M)^{\ast}_0$ denote the restriction maps ($\pi^i_j(h) = \pi^i(h) = h|_{M_i}$). 
Let ${\cal V}^i_j = (\pi^i)^{-1}({\cal U}^i_j)$. 

By Proposition 2.4 the maps $\pi^i_j$ and $\pi^i$ are principal bundles. 
Our stratigy is to use these bundles in order to reduce the noncompact case to the compact case. In \cite{Ya3} we have shown 

\begin{lemma} 
(i) ${\cal H}_j \cong {\cal H}_{X \cup L_j}(M_j)_0$ is an AR . \\
(ii) ${\cal U}^i_j$ is open in ${\cal E}_X(M_i, M)^{\ast}_0$, $cl \, {\cal U}^i_j \subset {\cal U}^i_{j+1}$ and ${\cal E}_X(M_i, M)^{\ast}_0 = \cup_{j > i} \, {\cal U}^i_j$. \\
(iii) ${\cal H}_X(M)_0 = \cup_{j > i}\,{\cal V}^i_j$, $cl \, {\cal V}^i_j \subset {\cal V}^i_{j+1}$ ($j > i$) and 
${\cal V}^{i+1}_j \subset {\cal V}^i_j$ ($j > i+1$). \\
(iv) The restriction map  $\pi^i_j : {\cal H}_j \to {\cal U}^i_j$ is a principal bundle with structure group ${\cal G}^i_j \equiv {\cal H}_j \cap {\cal H}_{M_i}(M) \cong {\cal H}_{X \cup L_j}(M_j)_0 \cap {\cal H}_{M_i \cup L_j}(M_j)$.  \\
(v) In the case (II), for each $j > i \geq 0$, (a) ${\cal G}^i_j$ is an AR $(\cong {\cal H}_{M_i \cup L_j}(M_j)_0)$, (b) $\pi^i_j : {\cal H}_j \to {\cal U}^i_j$ is a trivial bundle and (c) ${\cal U}^i_j$ is also an AR. 
\end{lemma}

\begin{Proof of Proposition 3.1}{\bf Case (II)}: \hspace{4pt} The assertion will be verified in the next four steps: \\
(1) For each $i \geq 1$ there exists a map $s^i : {\cal E}_X(M_i, M)^{\ast}_0 \to {\cal H}_X^{\, c}(M)_0$ such that $s^i(h)|_{M_i} = h \ (h \in {\cal E}_X(M_i, M)^{\ast}_0)$ and $s^i(cl \, {\cal U}^i_j) \subset {\cal H}_{j+1} \ (j > i)$. \\
(2) There exists a homotopy $G : {\cal H}_X(M)_0 \times [0,1] \to {\cal H}_X(M)_0$ such that (a) $G_0 = id$ and $G_t({\cal H}_X(M)_0)$ $\subset {\cal H}_X^{\, c}(M)_0 \ (0 < t \leq 1)$ and (b) each $(h, t) \in {\cal H}_X(M)_0 \times (0,1]$ admits a neighborhood ${\cal W}$ such that $G({\cal W}) \subset {\cal H}_i$ for some $i \geq 1$. \\
(3) There exists a function $\Phi : {\cal H}_X^{\, c}(M) \times [0, 1] \to {\cal H}_X^{\, c}(M)$ such that $\Phi_0 = id$, $\Phi_t({\cal H}_X^{\, c}(M)) \subset {\cal H}^{{\rm PL},c}_X(M) \ (0 < t \leq 1)$ and $\Phi|_{{\cal H}_{X \cup (M \setminus U_i)}(M) \times [0, 1]}$ is continuous for each $i \geq 1$. \\
(4) There exists a homotopy $F : {\cal H}_X(M)_0 \times [0, 1] \to {\cal H}_X(M)_0$ such that $F_0 = id$ and $F_t({\cal H}_X(M)_0) \subset {\cal H}^{{\rm PL, c}}_X(M)_0 \ (0 < t \leq 1)$. 
\end{Proof of Proposition 3.1}

\begin{proof}
(1) Since the restriction map $\pi^i_j : {\cal H}_j \to {\cal U}_j^i$ is a trivial bundle with an AR fiber and ${\cal H}_j \subset {\cal H}_{j+1}$, inductively we can construct a section $s^i_j : cl \, {\cal U}^i_j \to {\cal H}_{j+1}$ of $\pi^i_{j+1}$ such that $s_{j+1}^i|_{cl \, {\cal U}^i_j} = s^i_j \ (j > i)$. The section $s^i$ is defined by $s^i|_{cl \, {\cal U}^i_j} = s^i_j$. 

(2) We replace the interval $[0,1]$ by $[1, \infty]$. For each $i \geq 1$ let $G_i = s^i \pi^i : {\cal H}_X(M)_0 \to {\cal H}_X^{\, c}(M)_0$. We have $G_i(cl \, {\cal V}^i_j) \subset {\cal H}_{j+1}$ and $G_i(h)|_{M_i} = h|_{M_i}$. Since $\pi^i_j : {\cal H}_j \to {\cal U}_j^i$ is a trivial bundle with an AR fiber, using the \lq\lq fiber preserving" absolute extension property, we can inductively construct a sequence of homotopies $g^j : cl \, {\cal V}_j^{i+1} \times [i, i+1] \to {\cal H}_{j+1} \ (j > i+1)$ such that $g^j_i = G_i$, $g^j_{i+1} = G_{i+1}$, $g^{j+1}|_{cl \, {\cal V}_j^{i+1} \times [i, i+1]} = g^j$ and $g^j_t(h)|_{M_i} = h|_{M_i}$. Hence we can define a homotopy $G : {\cal H}_X(M)_0 \times [i, i+1] \to {\cal H}_X^{\, c}(M)_0$ by $G|_{cl \, {\cal V}_j^{i+1} \times [i, i+1]} = g^j$. Since $G_t(h)|_{M_i} = h|_{M_i}$ for $t \geq i$, we can continuously extend G by $G_{\infty} = id$. Since $G(cl \, {\cal V}_j^{i+1} \times [i, i+1]) \subset {\cal H}_{j+1} \ (j > i+1)$, the homotopy satisfies the required conditions. 

(3) Since $(X_i, Y_i) \equiv ({\cal H}_{X \cup (M \setminus U_i)}(M), {\cal H}_{X \cup (M \setminus U_i)}^{{\rm PL}}(M)) \cong ({\cal H}_{X \cup L_i}(M_i), {\cal H}_{X \cup L_i}^{{\rm PL}}(M_i))$, 
from Proposition 2.2 and Proposition 3.1 (Compact case) it follows that $X_i$ is an ANR and $Y_i$ has the h.n.\ complement in $X_i$. Inductively we can construct a sequence of homotopies $\phi^i : X_i \times [0, 1] \to X_i \ (i \geq 1)$ such that $\phi^i_0 = id$, $\phi^i_t(X_i) \subset Y_i \ (0 < t \leq 1)$ and $\phi^{i+1}|_{X_i \times [0, 1]} = \phi^i$. 
In fact, since $X_{i+1}$ is an ANR, given $\phi^i$, we have a homotpy $\psi_t : X_{i+1} \to X_{i+1}$ such that $\psi_0 = id$ and $\psi|_{X_i \times [0, 1]} = \phi^i$. 
Let $f_t$ be an absorbing homotopy of $X_{i+1}$ into $Y_{i+1}$ and take a map $\alpha : X_{i+1} \times [0, 1] \to [0, 1]$ with $\alpha^{-1}(0) = X_{i+1} \times \{0\} \cup X_i \times [0, 1]$. 
The desired homotopy $\phi^{i+1}$ is defined by $\phi^{i+1}(x, t) = f(\psi(x, t), \alpha(x, t))$. Since ${\cal H}_X^{\, c}(M) = \cup_i X_i$ and ${\cal H}_X^{{\rm PL},c}(M) = \cup_i Y_i$, the required function $\Phi$ is defined by $\Phi|_{X_i \times [0, 1]} = \phi^i \ (i \geq 1)$. 

(4) We define a function $H : {\cal H}_X(M)_0 \times (0, 1] \times [0, 1] \to {\cal H}_X^{\, c}(M)$ by $H(h, t, s) = \Phi(G(h, t), s)$. 
By (2) and (3) it follows that $H$ is continuous and $H(h, t, 0) = G(h, t)$. 
Choose a small map $\alpha : {\cal H}_X(M)_0 \times (0, 1] \to (0, 1]$ such that 
$\rho(H(h,t,\alpha(h,t)), G(h,t)) < t$ ($(h, t) \in {\cal H}_X(M)_0 \times (0, 1]$) 
and define the homotopy $F : {\cal H}_X(M)_0 \times [0, 1] \to {\cal H}_X(M)$ by 
\[ F(h, t) = \begin{cases}
\ \ h \ \ (= G_0(h)) &  \ t = 0 \\
H(h, t, \alpha(h, t)) & 0 < t \leq 1 
\end{cases} \]
Then $F$ is continuous and it follows that ${\rm Im} \, F \subset {\cal H}_X(M)_0$ and $F({\cal H}_X^{{\rm PL}, c}(M)_0 \times [0, 1] \cup {\cal H}_X(M)_0 \times (0, 1]) = {\cal H}_X^{{\rm PL}, c}(M)_0$. 
This completes the proof.
\vskip 2mm
{\bf Case (I)}: In this case ${\cal H}_X(M)$ is an ANR, and the conclusion follows from the next lemma. 
\end{proof}

\begin{lemma}
If ${\cal H}_X(M)$ is an ANR, then ${\cal H}_X^{{\rm PL, c}}(M)_0$ has the h.n.\ complement in ${\cal H}_X(M)_0$.
\end{lemma} 

\begin{proof}
We verify the following assertion: 

\renewcommand{\labelitemi}{$(\ast)$}
\begin{itemize}
\item Every $h \in {\cal H}_X(M)_0$ admits an open neighborhood ${\cal V}$ and 
a homotopy $\phi : {\cal V} \times [0, 1] \to {\cal H}_X(M)_0$ such that $\phi_0$ 
is the inclusion ${\cal V} \subset {\cal H}_X(M)_0$ and $\phi_t({\cal V}) \subset {\cal H}_X^{{\rm PL, c}}(M)_0$ $(0 < t \leq 1)$.  
\end{itemize}

\noindent This easily implies that ${\cal H}_X(M)_0$ itself admits an absorption homotopy into ${\cal H}_X^{{\rm PL, c}}(M)_0$ (cf. \cite[Fact 4.1 (i)]{Ya3}). 
Consider the bundle ${\cal H}_{M_1}(M) \subset {\cal H}_X(M) \to {\cal E}_X(M_1, M)^{\ast}$ (over its image) 
(cf. Proof of \cite[Corollary 1.1]{Ya3}, \cite[Corollary 2.1]{Ya2}). 
From the assumption it follows that ${\cal H}_{M_1}(M)$ is also an ANR, in which ${\cal H}_{M_1}(M)_0$ is open, and ${\cal H}_{M_i}(M) \subset {\cal H}_{M_1}(M)_0$ for some $i \geq 1$ since ${\rm diam} \, {\cal H}_{M_i}(M) \leq 1/2^i$. 
Let $h \in {\cal H}_X(M)_0$ and consider the restriction map $\pi^i : {\cal H}_X(M)_0 \to {\cal E}_X(M_i, M)^\ast_0$. 
By Lemma 3.1 (ii) $h|_{M_i} \in {\cal U}_j^i$ for some $j > i$. 
Since the map $\pi^i_j : {\cal H}_j \to {\cal U}^i_j$ is locally trivial, it has a section $s$ over an open neighborhood ${\cal U}$ of $h|_{M_i}$ in ${\cal U}_j^i$. Then ${\cal V} = (\pi^i)^{-1}({\cal U})$ is an open neighborhood of $h$ in ${\cal H}_X(M)_0$ and for any $g \in {\cal V}$ we have $s(g|_{M_i})^{-1} g \in {\cal H}_{M_i}(M) \subset {\cal H}_{M_1}(M)_0$. 
By Proposition 3.1 Case (II), ${\cal H}_{M_1}^{{\rm PL, c}}(M)_0$ has the h.n.\ complement in ${\cal H}_{M_1}(M)_0$. 
Since $({\cal H}_j, {\cal H}_{X \cup (M \setminus U_j)}^{{\rm PL}}(M)_0) \cong ({\cal H}_{X \cup L_j}(M_j)_0,{\cal H}_{X \cup L_j}^{{\rm PL}}(M_j)_0)$, 
by Proposition 3.1 (Compact case), ${\cal H}_{X \cup (M \setminus U_j)}^{{\rm PL}}(M)_0$ also has the h.n.\ complement in ${\cal H}_j$. 
Let $\psi_t$ and $\chi_t$ are the corresponding absorbing homotopies. 
Then the required absorbing homotopy $\phi_t$ is defined by $\phi_t(g) = \chi_t(s(g|_{M_i})) \psi_t(s(g|_{M_i})^{-1}g) \ (g \in {\cal V})$. 
\end{proof} 

\begin{example}
(1) In general, ${\cal H}_X^{\, c}(M) = \cup_{i \geq 1}{\cal H}_{X \cup (M \setminus U_i)}(M)$, but ${\cal H}_X^{\, c}(M)_0 
\supsetneqq 
\cup_{i \geq 1}{\cal H}_i$. For example, let $(M, X) = ({\Bbb S}^1 \times [0, \infty), {\Bbb S}^1 \times \{0\})$ and define $h_t \in {\cal H}_X^{\, c}(M) \ (1 \leq t < \infty)$ as the Dehn twist on the annulus ${\Bbb S}^1 \times [t, t+1]$ and $h_{\infty} = id_M$. Then $h_t \ (1 \leq t \leq \infty)$ is a path in ${\cal H}_X^{\, c}(M)_0$. However one can easily see that $h_t \not\in \cup_{i \geq 1}{\cal H}_i$ for any $1 \leq t < \infty$. 

(2) We will identify ${\Bbb R}^2$ with ${\Bbb C}$. The rotation map $f : {\Bbb S}^1 \to {\cal H}({\Bbb R}^2)_0$, $f(z)(x) = zx \ (z \in {\Bbb S}^1, x \in {\Bbb C})$ is a homotopy equivalence (cf \cite{Ha}). By Proposition 3.1 we can absorb $f$ into ${\cal H}^c({\Bbb R}^2)_0 (= {\cal H}^c({\Bbb R}^2))$. In fact, we can explicitly define an absorbing homotopy $F : {\Bbb S}^1 \times [0, 1] \to {\cal H}({\Bbb R}^2)_0$ such that $F_0 = f$, $F_t({\Bbb S}^1) \subset {\cal H}^c({\Bbb R}^2)_0 \ (0 < t \leq 1)$ in the following way: Choose any homeomorphism $\alpha : [0, 2\pi) \cong [1, \infty)$ and, replacing $[0, 1]$ by $[0, \infty]$, define $F$ by (a) $F_{\infty} = f$ and (b) for $z = e^{i \theta} \ (0 \leq \theta < 2\pi)$, $0 \leq t < \infty$ and $x \in {\Bbb R}^2$ 
\[ F(z, t)(x) = \begin{cases}
zx & \hspace{10pt}|x| \leq \alpha(\theta) + t, \\
e^{i s \theta}x & \hspace{10pt} \alpha(\theta) + t \leq |x| \leq \alpha(\theta) + t + 1 \ \mbox{and} \ s = \alpha(\theta) + t + 1 - |x|, \\
x & \hspace{10pt} |x| \geq \alpha(\theta) + t + 1.
\end{cases} \] 
\end{example}

\section{Stability and class property of homeomorphism groups} 
\subsection{Class property of homeomorphism groups} \mbox{} 

When $(X, d)$ and $(Y, \rho)$ are metric spaces, a map $f : X \to Y$ is said to be $K$-Lipschitz ($K \geq 1$) if $\rho(f(x), f(y)) \leq K d(x, y)$ for any $x, y \in X$. A map $f$ is Lipschitz if it is $K$-Lipschitz for some $K \geq 1$, and $f$ is locally Lipschitz if each $x \in X$ has a neighborhood $U$ such that $f|_U : U \to Y$ is Lipschitz. Every locally Lipschitz map is Lipschitz over any compact subsets. A (locally) Lipschitz homeomorphism is a homeomorphism $f$ such that both $f$ and $f^{-1}$ are (locally) Lipschitz. A Lipschitz embedding is a Lipschitz homeomorphism onto its image. 

A metric space is proper if any closed bounded subset is compact. A Euclidean polyhedron $X$ is a (closed) subpolyhderon of a Euclidean space with the induced metric. It follows that $X$ is a proper metric space and that ${\cal H}^{\rm PL}(X) \subset {\cal H}^{\rm loc \, LIP}(X)$ \cite{SW}.

\begin{lemma} 
(1) If $X$ is a locally compact polyhedron and $A$ is a closed subset of $X$, 
then ${\mathcal H}^{\rm PL}_A(X)$ is $F_{\sigma\delta}$ in ${\mathcal H}_A(X)$ and ${\mathcal H}^{\rm PL,c}_A(X)$ is $\sigma$-fd-compact 
\cite{Ge2}, \cite[Lemma 3.14]{Ya1}. \\ 
(2) Suppose $X = (X, d)$ is a locally compact, locally connected separable metric space and $A$ is a closed subset of $X$. 
(i) ${\mathcal H}^{\rm loc \, LIP}_A(X)$ is $F_{\sigma\delta}$ in ${\mathcal H}_A(X)$, 
and ${\mathcal H}^{\rm LIP, c}_A(X)$ is $\sigma$-compact \cite[Lemma 3.14]{Ya1}. \\
(ii) If $X$ is a proper metric space, then ${\mathcal H}^{\rm LIP}_A(X)$ is also $\sigma$-compact. 
\end{lemma}

\begin{proof}
(2)(ii) Choose a point $x_0 \in X$ and let $B_r = \{ x \in X \mid d(x, x_0) \leq r\}$ ($r \geq 0$). 
Since $X$ is proper, $B_r$ is compact.
For $K \geq 1$ and $r \geq 0$ consider the subspace 
\[ {\cal H}^K_r = \left\{ f \in {\cal H}^{\rm LIP}(M) \mid \frac{1}{K}d(x, y) \leq d(f(x), f(y)) \leq Kd(x, y) \, (x, y \in X), \, f(x_0), f^{-1}(x_0) \in B_r \right\}. \]
Since $\displaystyle {\cal H}^{\rm LIP}(M) = \cup_{n = 1}^\infty\,{\cal H}^n_n$, 
it suffices to show that each ${\cal H}^K_r$ is compact.
For each $f \in {\cal H}^K_r$ we have $f(B_k)$, $f^{-1}(B_k) \subset B_{r_k}$, where $r_k = r + Kk$.
By the Ascoli-Arzel\`a theorem, the space of $K$-Lipschitz maps from $B_k$ to $B_{r_k}$, ${\cal C}^{K \text{-LIP}}(B_k, B_{r_k})$, is compact, and the map
\[ \phi : {\cal H}^K_r \to \left( \prod_{k=1}^\infty {\cal C}^{K \text{-LIP}}(B_k, B_{r_k}) \right)^2, \ \phi(f) = ((f|_{B_k})_k, (f^{-1}|_{B_k})_k) \]
is a closed embedding. Thus ${\cal H}^K_r$ is compact as required.
\end{proof}

Suppose $M$ is a connected Riemannian $n$-manifold. 
We always assume that $M$ is equipped with the path metric $d_M$ induced from the Riemannian metric of $M$.
We say that $M$ is complete if the metric $d_M$ is complete. 
In this case $(M, d_M)$ is a proper metric space \cite[Hopf - Rinow Theorem p.~111]{Sak}. 
Suppose $M$ has a $C^1$-triangulation. 
As a polyhedron, $M$ admits a PL-homeomorphism $h$ onto a connected Euclidean polyhedron $Y \subset {\Bbb R}^N$ ($N = 2n+1$). 
The standard metric of ${\Bbb R}^N$ is denoted by $d$.

\begin{lemma}
The PL-homeomorphism $h : (M, d_M) \to (Y, d)$ is a locally Lipschitz homeomorphism.
\end{lemma}

\begin{proof}
We define a compatible metric $\rho_Y$ of $Y$ by 
\[ \mbox{$\rho_Y(x, y) = \inf\,\{ L(\alpha) \mid \alpha$ is a PL-path in $Y$ joining $x$ and $y\}$} \]
where $L(\alpha)$ is the length of a path $\alpha$ in $({\Bbb R}^N, d)$. 

First we note that $id : (Y, d) \to (Y, \rho_Y)$ is a locally Lipschitz homeomorphism. 
This follows from the following observations: 
(i) We may assume that $Y$ is compact. 
(ii) In general, if $K$ and $L$ are compact connected Euclidean polyhedra and $f : K \to L$ is a PL-map, 
then $f : (K, \rho_K) \to (L, \rho_L)$ is a Lipschitz map. 
(iii) Take any regular neighborhood $N$ of $Y$ in ${\Bbb R}^N$ and a PL-retraction $r : N \to Y$.  
Then $N$ contains a $\varepsilon$-neighborhood $N(Y, \varepsilon)$ for some small $\varepsilon > 0$. 
It follows that $r : (N, \rho_N) \to (Y, \rho_Y)$ is Lipschitz and $\rho_N(x, y) = d(x, y)$ for $x, y \in Y$ with $d(x, y) < \varepsilon$, 
so that $id : (Y, d) \to (Y, \rho_Y)$ is also Lipschitz.

It remains to show that $h : (M, d_M) \to (Y, \rho_Y)$ is a locally Lipschitz homeomorphism. 
Let $T$ denote the $C^1$-triangulation of $M$. Subdividing $T$ if necessary, we may assume that $h$ is simplicial on each simplex of $T$. 
Since $T$ is a $C^1$-triangulation, for each $n$-simplex $\sigma \in T$ the restriction $h : \sigma \to h(\sigma)$ is $C^1$-diffeomorphism 
and there is a constant $K_\sigma \geq 1$ such that for each piecewise $C^1$-curve $\alpha$ in $\sigma$
\[ \frac{1}{K_\sigma} L_M(\alpha) \leq L(h(\alpha)) \leq K_\sigma L_M(\alpha). \]

We show that $h$ is locally Lipschitz. A similar argument shows that $h^{-1}$ is also locally Lipschitz. 
Let $x$ be any point of $M$ and $U$ be a strongly convex open neighborhood of $x$ in $M$ (cf. \cite[Ch. IV, \S 5]{Sak}) containd in a compact subcomplex $S$ of $T$. 
Let $K = \max\{ K_\sigma \mid \sigma^n \in S\}$. 
Any points $y$ and $z$ in $U$ can be joined by a unique shortest geodesic $\alpha$ of $M$ contained in $U$. 
For any $\varepsilon > 0$ we can find a piecewise $C^1$-curve $\beta$ in $U$ joining $y$ and $z$ such that $L_M(\beta) < L_M(\alpha) + \varepsilon$ and $\beta$ decomposes into subcurves $\beta_i$ so that each $\beta_i$ is contained in an $n$-simplex $\sigma_i \in S$ (cf. \cite[Ch.3. Theorem 2.1]{Hir}). 
Then $L(h(\beta)) \leq K\, L_M(\beta)$. 
Since $d_M(y, z) = L_M(\alpha)$ and $\rho_Y(h(y), h(z)) \leq L(h(\beta))$, it follows that $\rho_Y(h(y), h(z)) \leq K\,d_M(y, z)$. 
\end{proof}

Suppose $X$ is a compact subpolyhedron of $M$ and let $A = h(X)$. By Lemma 4.2 the homeomorphism $h$ induces a homeomorphism of tuples 
\begin{eqnarray*}
h^{\#} & : & ({\mathcal H}_X(M), {\mathcal H}^{\rm loc\,LIP}_X(M), {\mathcal H}^{\rm LIP, c}_X(M), 
{\mathcal H}^{\rm PL}_X(M), {\mathcal H}^{\rm PL, c}_X(M)) \\
& & \hspace{20mm} \stackrel{\cong}{\longrightarrow} \hspace{5mm} ({\mathcal H}_A(Y), {\mathcal H}^{\rm loc\,LIP}_A(Y), {\mathcal H}^{\rm LIP, c}_A(Y), 
{\mathcal H}^{\rm PL}_A(Y), {\mathcal H}^{\rm PL, c}_A(Y))
\end{eqnarray*}

If $h$ is a Lipschitz homeomorphism, then we can include ${\mathcal H}^{\rm LIP}_X(M)$ in the above tuple. 
However, the author does not know whether every complete Riemannian manifold with a $C^1$-triangulation 
admits a Lipschitz PL-embedding into a Euclidean space (or a Hilbert space). 

\subsection{Stability property of homeomorphism groups} \mbox{}

K.\,Sakai and R.Y.\,Wong showed that the triple $({\cal H}(X), {\cal H}^{{\rm LIP}}(X) {\cal H}^{\rm PL}(X))$ is 
$(s, \Sigma, \sigma)$-stable for every Euclidean polyhedron $X$, by using the Morse's $\mu$-length of arcs \cite{SW}. 
In the case where $Y$ is noncompact we can verify a more precise statement. 
In this subsection we do {\it not} assume that a tuple $(X, X_1, \cdots, X_\ell)$ satisfies the inclusion relation: $X_i \supset X_j$ ($1 \leq i < j \leq \ell$)

\begin{proposition}
Suppose $X$ is a Euclidean polyhedron and $A$ is a compact subpolyhedron of $X$ such that ${\rm dim}\,(X \setminus A) \geq 1$. \\
(i) $\displaystyle ({\mathcal H}_A(X), {\mathcal H}^{\rm loc\,LIP}_A(X), {\mathcal H}^{\rm LIP}_A(X), {\mathcal H}^{\rm LIP, c}_A(X), 
{\mathcal H}^{\rm PL}_A(X), {\mathcal H}^{\rm PL, c}_A(X))$ \ is $(s, \Sigma, \Sigma, \Sigma, \sigma, \sigma)$-stable. \\ 
(ii) If $X$ is noncompact and ${\rm dim}\,(X \setminus K) \geq 1$ for any compact subset K of X 
(or equivalently, any triangulation of $X$ contains infinitely many 1-simplices), then this tuple is 
$(s^{\infty}, \Sigma^{\infty}, s^\infty_b, \Sigma^{\infty}_f, \sigma^{\infty}, \sigma^{\infty}_f)$-stable 
\end{proposition}

For Riemannian manifolds with $C^1$-triangulations we have the following version:

\begin{proposition}
Suppose $M$ is a connected Riemannian $n$-manifold with a $C^1$-triangulation and $X$ is a compact subpolyhedron of $M$ with $X \neq M$. \\
(i) $({\mathcal H}_X(M), {\mathcal H}^{\rm loc \, LIP}_X(M), {\mathcal H}^{\rm LIP}_X(M), 
{\mathcal H}^{\rm LIP, c}_X(M), {\mathcal H}^{\rm PL}_X(M), {\mathcal H}^{\rm PL, c}_X(M))$ is $(s, \Sigma, \Sigma, \Sigma, \sigma, \sigma)$-stable. \\ 
(ii) If $M$ is noncompact, then $({\mathcal H}_X(M), {\mathcal H}^{\rm loc\,LIP}_X(M), {\mathcal H}^{\rm LIP, c}_X(M), {\mathcal H}^{\rm PL}_X(M), {\mathcal H}^{\rm PL, c}_X(M))$ is 
$(s^{\infty}, \Sigma^{\infty}, \Sigma^{\infty}_f$, $\sigma^{\infty}, \sigma^{\infty}_f)$-stable. \\
(iii) If $M$ is noncompact and contains a sequence of points $p_k$ ($k = 1, 2, \cdots$) 
such that $d_M(p_k, p_\ell) \geq \varepsilon$ ($k \neq \ell$) for some constant $\varepsilon > 0$, then $({\cal H}_X(M), {\cal H}_X^{\rm loc \, LIP}(M), {\cal H}_X^{\rm LIP}(M),{\cal H}_X^{\rm LIP, c}(M))$ is $(s^{\infty},\Sigma^{\infty}, s^\infty_b, \Sigma^{\infty}_f)$-stable. 
\end{proposition}

In Propositions 4.1 and 4.2 each stability homeomorphism preserves the identity component part. 

The statements (i), (ii) in Proposition 4.2 follows from Proposition 4.1 (i), (ii) and Lemma 4.2, 
while the statement (iii) follows from the direct argument similar to the proof of Proposition 4.1,
with replacing straight segments by geodesics. This direct argument does not enable us to include the subgroups of PL-homeomorphisms. 

In order to verify Proposition 4.1 and 4.2 in a systematic way, we have to clarify the stability argument in \cite{SW} and extract some conditions under which the stability homeomorphisms preserve the subspaces of tuples in Proposition 4.1 and 4.2.  
For this purpose first we recall the definition and basic properties of the Morse's $\mu$-length of arcs: 

Suppose $(X, d)$ is a metric space and $A$ is an arc (or a point) in $X$. 
The arc $A$ admits a canonical linear order $\leq$ unique up to the inversion. 
For each $k \geq 1$ let $S_k = \{ (a_0, a_1, \cdots, a_k) \in A^{k+1} \mid a_0 \leq a_1 \leq \cdots \leq a_k\}$ and define
$\delta({\boldsymbol a}) = \min\{d(a_{i-1}, a_i) \mid i = 1, \cdots, k\}$ (${\boldsymbol a} = (a_0, \cdots, a_k) \in S_k$) and 
$\mu_k(A) = \sup\{\delta({\boldsymbol a}) \mid {\boldsymbol a} \in S_k\}$.
The $\mu$-length of $A$ is defined as $\mu(A) = \sum_{k=1}^\infty 2^{-k} \mu_k(A)$. If $A$ is a point, then $\mu(A) = 0$.
For example, $\displaystyle \mu([a, b]) = (b-a) \sum_{n=1}^\infty \frac{1}{n \, 2^n}$ for any closed interval $[a, b] \subset {\Bbb R}$.
The following facts are easily verified from the definition (\cite[\S1.~pp.~197--202]{SW}).

\begin{lemma}
(i) If $A_1$ is a subarc of $A$ and $A_1 \neq A$, then $\mu(A_1) < \mu(A)$. \\
(ii) The functions ${\cal E}([-1, 1], X) \times [-1, 1] \to [0, \infty)$, $(f, t) \mapsto \mu(f([-1, t]))$, $(f, t) \mapsto \mu(f([t, 1]))$, 
are continuous. \\
(iii) For each $f \in {\cal E}([-1, 1], X)$, there is a unique $t_f \in (-1, 1)$ such that $\mu(f([-1, t_f])) = \mu(f([t_f, 1]))$. 
The function ${\cal E}([-1, 1], X) \to (-1, 1)$, $f \mapsto t_f$, is continuous. \\
(iv) Suppose $f : (A, d_A) \to (X, d_X)$ is an embedding of an arc $A$. \\
(a) If $f$ is $K$-Lipschitz (i.e., $K > 0$ and $d_Y(f(x), f(y) \leq K\,d_A(x,y)$ ($x, y \in X$)), then $\mu(f(A)) \leq K\,\mu(A)$. \\
(b) If $f$ is $c$-similar (i.e., $c > 0$ and $d_Y(f(x), f(y)) = c \, d_X(x, y)$ ($x, y \in X$)), then $\mu(f(A)) = c\,\mu(A)$.\\
(v) If $f \in {\cal E}([-1, 1], X)$ is a $c$-similar embedding ($c > 0$) and $g \in {\cal E}(X, Y)$ is a $K$-Lipschitz embedding ($K \geq 1$), 
then $t_f = 0$ and $\displaystyle |t_{gf}| \leq \frac{K^2 - 1}{K^2 + 1} < 1$.
\end{lemma}

\begin{proof}
(v) The value $t = t_{gf}$ satisfies the following equations: 
\begin{alignat*}{2}
& \hspace{5mm} \mu(gf([-1, t])) = \mu(gf([t, 1])), & & \text{(Definition of $t$)} \\[2mm]
& \mbox{\small 
$\displaystyle 
\left\{\begin{array}{l}
\mbox{$\displaystyle \frac{1}{K}$} \mu(f([-1, t])) \leq \mu(gf([-1, t])) \leq K \mu(f([-1, t])) \\[4mm]
\mbox{$\displaystyle \frac{1}{K}$} \mu(f([t, 1])) \leq \mu(gf([t, 1])) \leq K \mu(f([t, 1]))
\end{array}
\right.
$ \hspace{5mm}}
& & \text{($g$ is a $K$-Lipschitz embedding)} \\[2mm]
& \mbox{\small 
$\displaystyle
\left\{\begin{array}{l}
\mu(f([-1, t])) = c (t + 1) \sum_{m=1}^\infty \frac{1}{m\,2^m} \\[4mm]
\mu(f([t, 1])) = c(1 - t) \sum_{m=1}^\infty \frac{1}{m\,2^m}
\end{array}
\right.
$}
& & \text{($f$ is $c$-similar)} 
\end{alignat*}
Thus {\small $\displaystyle \frac{1}{K} (t+1) \leq K (1 - t)$} and {\small $\displaystyle \frac{1}{K} (1 - t) \leq K (t+1)$}, 
and so {\small $\displaystyle |t| \leq \frac{K^2 - 1}{K^2 + 1}$}.
\end{proof}

\begin{example}
(i) If $X$ is a normed linear space and $f \in {\cal E}([a, b], X)$ is affine linear, then $f$ is $c$-similar for $c = \|f(b) - f(a)\|/(b-a)$. \\
(ii) If $M$ is a connected Riemannian manifold and $f \in {\cal E}([a, b], M)$ is a shortest geodesic of constant velocity $c > 0$, 
then $f$ is a $c$-similar embedding.
\end{example}

Next we recall the construction of the stability homeomorphism in \cite{SW}. 
Below we regard as $s = \prod_{i=1}^\infty(-1, 1)$ (instead of ${\Bbb R}^\infty$) 
so that $\Sigma = \{(x_n) \in s \, : \, \sup_n\,|x_n| < 1 \}$. 

In ${\Bbb R}^n$ we use the norm defined by $\|x\|_1 = \sum_{i = 1}^n |x_i|$ ($x = (x_i) \in {\Bbb R}^n$). 
Let $B^n = \{ x \in {\Bbb R}^n \mid \| x \|_1 \leq 1 \}$. 
We identify the interval $[-1, 1]$ with the segment in $B^n$ connecting the two vertices $(\pm 1, 0,\cdots, 0)$. 
For $t \in (-1, 1)$ we define $\overline{\lambda}(t) \in {\cal H}_\partial^{\rm PL}(B^n)$ by 
\[ \overline{\lambda}(t)(x) = x + ((1 - \| x \|_1)~t, 0, \cdots, 0) \hspace{10mm}  (x \in B^n). \]
The function $\overline{\lambda} : (-1, 1) \to {\cal H}_\partial^{\rm PL}(B^n)$ is continuous. 

For each $i \geq 1$ 
let $a_i = (\frac{1}{2^i}, 0 ,\cdots, 0) \in {\Bbb R}^n$ and $c_i = \frac{1}{2^{i+2}}$, and  
define a $c_i$-similar transformation $\alpha_i : {\Bbb R}^n \to {\Bbb R}^n$ by $\alpha_i(x) = c_i x + a_i$. 
Let $B_i^n = \alpha_i(B^n)$. 
Note that 
$B_i^n = \{ x \in {\Bbb R}^n \mid \| x - a_i \|_1 \leq c_i \} \subset {\rm Int}\,B^n$, 
$\alpha_i([-1, 1]) = B_i^n \cap [-1, 1]$ and $B_i^n \cap B_j^n = \emptyset$ for $i \neq j$. 
For $t \in (-1, 1)$ we define $\lambda_i(t) \in {\cal H}_\partial^{\rm PL}(B_i^n)$ by $\lambda_i(t) = \alpha_i \overline{\lambda}(t) \alpha_i^{-1}$.

For each $\mbox{\boldmath $t$} = (t_i) \in s$, 
we define $\lambda({\boldsymbol t}) \in {\cal H}_\partial(B^n)$ by $\lambda({\boldsymbol t}) = \lambda_i(t_i)$ on $B_i^n$ 
and $\lambda({\boldsymbol t}) = id$ on $B^n \setminus \cup_{i = 1}^\infty B_i^n$. 
The function $\lambda : s \to {\cal H}_\partial(B^n)$ is continuous. \\
(I) Suppose $X$ is a space, $Y$ is a metric space, 
$\phi : B^n \to X$ is an embedding such that $\phi({\rm Int}\, B^n)$ is open in $X$, and $\chi : s^2 \cong s$ is a homeomorphism.  
Let ${\cal F} = \{ f \in {\cal C}(X, Y) \mid f : \phi([-1, 1]) \to Y$ is an embedding$\}$. 
For these data, the stability homeomorphism $\Phi : {\cal F} \times s \to {\cal F}$ is defined as follows: \\ 
(1) (a) For each $f \in {\cal F}$ and $i \geq 1$ there exists a unique $t = t_i(f) \in (-1, 1)$ 
such that $\mu(f(\phi(\alpha_i([-1,t])))) = \mu(f(\phi(\alpha_i([t,1]))))$. 
Define a map $F : {\cal F} \to s$ by $F(f) = (t_i(f))_i$. \\
(b) Define a map $G : s \to {\cal H}_{X \setminus \phi({\rm Int}\,B^n)}(X)$ by 
$G({\boldsymbol t}) = \phi \lambda({\boldsymbol t}) \phi^{-1}$ on $\phi(B^n)$ and 
$G(t) = id$ on $X \setminus \phi({\rm Int}\,B^n)$ for ${\boldsymbol t} \in s$. \\
(2) We have reciplocal homeomorphisms 
\[ \begin{cases} 
\Phi' : {\cal F} \to F^{-1}({\boldsymbol 0}) \times s, & \hspace{5mm} \Phi'(f) = (f G(F(f)), F(f)) \\ 
\Psi' : F^{-1}({\boldsymbol 0}) \times s \to {\cal F}, & \hspace{5mm} \Psi'(f, {\boldsymbol t}) = f G({\boldsymbol t})^{-1} 
\end{cases}
\]
The stability homeomorphism $\Phi : {\cal F} \times s \to {\cal F}$ is defined as the composition 
\[ {\cal F} \times s \stackrel{\Phi' \times 1}{\longrightarrow} (F^{-1}({\boldsymbol 0}) \times s) \times s \cong 
F^{-1}({\boldsymbol 0}) \times s^2 \stackrel{1 \times \chi}{\longrightarrow} F^{-1}({\boldsymbol 0}) \times s 
\stackrel{\Psi'}{\longrightarrow} {\cal F} \]

To show the stability property in the noncompact case we use a sequence of embeddings of $B^n$ into $X$. \\
(II) Suppose $X$ is a noncompact space, $Y$ is a metric space,  
$\phi_k : B^{n_k} \to X$, $n_k \geq 1$, ($k \geq 1$) is a sequence of embeddings such that 
$\phi_k({\rm Int}\, B^{n_k})$ is open in $X$, $\phi_k(B^{n_k}) \cap \phi_l(B^{n_l}) = \emptyset$ ($k \neq l$) and 
$\{ \phi_k(B^{n_k}) \}_k$ is discrete in $X$,  
and $\chi : (s^\infty)^2 \cong s^\infty$ is a homeomorphism. 
Let ${\cal F} = \{ f \in {\cal C}(X, Y) \mid f : \phi_k([-1, 1]) \to Y$ 
is an embedding for each $k \geq 1\}$. 
The stability homeomorphism $\Phi : {\cal F} \times s^\infty \to {\cal F}$ is defined as follows: \\
(0) For each $k \geq 1$ and $i \geq 1$, as in the last paragraph 
we have the $c_i$-similar transformation $\alpha_{k,i} : {\Bbb R}^{n_k} \to {\Bbb R}^{n_k}$ and 
the maps $\overline{\lambda}_k : (-1, 1) \to {\cal H}_\partial^{\rm PL}(B^{n_k})$, 
$\lambda_{k, i} : (-1, 1) \to {\cal H}_\partial^{\rm PL}(B_i^{n_k})$ and 
$\lambda_k : s \to {\cal H}_\partial(B^{n_k})$. \\
(1) We define two maps $F : {\cal F} \to s^\infty$ and $G : s^\infty \to {\cal H}_{X \setminus \cup_k\,\phi_k({\rm Int}\,B^{n_k})}(X)$ by: \\
(i) $F(f) = ((t_{k,i})_i)_k \in s^\infty$, where $t_{k,i} = t_{f\phi_k\alpha_{k,i}}$  
(by definition, $t = t_{f\phi_k\alpha_{k,i}}$ is the unique $t \in (-1, 1)$ such that $\mu(f(\phi_k(\alpha_{k,i}([-1,t])))) = \mu(f(\phi_k(\alpha_{k,i}([t,1]))))$. \\
(ii) $G({\boldsymbol t}) = \phi_k \lambda_k({\boldsymbol t}_k) \phi_k^{-1}$ on $\phi_k(B^{n_k})$ and 
$G({\boldsymbol t}) = id$ on $X \setminus \cup_k \phi_k(B^{n_k})$ for ${\boldsymbol t} = ({\boldsymbol t}_k)_k \in s^\infty$. \\
(2) The reciplocal homeomorphisms $\Phi' : {\cal F} \to F^{-1}({\boldsymbol 0}) \times s^\infty$ 
and $\Psi' : F^{-1}({\boldsymbol 0}) \times s^\infty \to {\cal F}$, and 
the stability homeomorphism $\Phi : {\cal F} \times s^\infty \to {\cal F}$ are defined by the same formulae. \\
(3) The maps $F$, $G$ and $\Phi$ have the following properties: 

\begin{lemma} (i) Suppose ${\cal F}_1 \subset {\cal F}$ and $S_1 \subset s^\infty$ are subsets 
and ${\cal G}_1 \subset {\cal H}(X)$ is a subgroup. 
If (a) $F({\cal F}_1) \subset S_1$, 
(b) $G(S_1) \subset {\cal G}_1$, 
(c) $fg \in {\cal F}_1$ for any $f \in {\cal F}_1$ and $g \in {\cal G}_1$ and 
(d) $\chi(S_1^2) = S_1$, then 
$\Phi({\cal F}_1 \times S_1) = {\cal F}_1$. \\
(ii) When $Y = X$, if (e) $F(id_X) = {\boldsymbol 0}$, $G({\boldsymbol 0}) = id_X$ and $\chi({\boldsymbol 0}, {\boldsymbol 0}) = {\boldsymbol 0}$, 
then $\Phi(id_X, {\boldsymbol 0}) = id_X$. 
\end{lemma} 

\begin{proof}
(i) The assumption implies that 
$\Phi'({\cal F}_1) = (F^{-1}({\boldsymbol 0}) \cap {\cal F}_1) \times S_1$ and 
$\Psi'((F^{-1}({\boldsymbol 0}) \cap {\cal F}_1) \times S_1) = {\cal F}_1$. \\
(ii) Set ${\cal F}_1 = {\cal G}_1 = \{ id_X \}$ and $S_1 = \{ {\boldsymbol 0} \}$.
\end{proof}

\begin{lemma} 
(i) Let $f \in {\cal F}$ and $F(f) = ((t_{k,i})_i)_k \in s^\infty$. \\
(a) If $\phi_{\ell} : [-1, 1] \to X$ is a similar embedding and $f : \phi_{\ell}([-1, 1]) \to Y$ is a $K$-Lipschitz embedding, then $|t_{\ell,i}| \leq \frac{K^2 - 1}{K^2 +1}$ ($i \geq 1$). \\
(b) If $f\phi_{\ell} : [0, a_j + c_j] \to Y$ is a similar embedding, then $t_{\ell,i} = 0$ ($i \geq j$). \\
(ii) Let ${\boldsymbol t} = ({\boldsymbol t}_k)_k \in s^\infty$. \\
(a) If $\phi_{\ell} : B^{n_{\ell}} \to X$ is a similar ($L$-Lipschitz) embedding and ${\boldsymbol t}_{\ell} \in [-s, s]^\infty$ ($0 \leq s < 1$), 
then $G({\boldsymbol t}) : \phi_{\ell}(B^{n_{\ell}}) \to \phi_{\ell}(B^{n_{\ell}})$ is a $1/(1 - s)$-Lipschitz ($L^2/(1 - s)$-Lipschitz) homeomorphism. \\
(b) If $X$ is a polyhedron, $\phi_{\ell}$ is a PL-embedding and ${\boldsymbol t}_{\ell} \in \sigma$, 
then $G({\boldsymbol t}) : \phi_{\ell}(B^{n_{\ell}}) \to \phi_{\ell}(B^{n_{\ell}})$ is a PL-homeomorphism. \\
(c) If ${\boldsymbol t}_{\ell} = {\boldsymbol 0}$, then $G({\boldsymbol t}) = id$ on $\phi_{\ell}(B^{n_{\ell}})$.
\end{lemma}

\begin{proof}
The statement (i) follows from Lemma 4.3~(v) and the definition of $t_{\ell,i}$. 
The statement (ii) follows from the next facts: \\
(a) If ${\boldsymbol t}_{\ell} \in [-s, s]^\infty$ ($0 \leq s < 1$), 
then $\lambda_{\ell}({\boldsymbol t}_{\ell}) : B^{n_{\ell}} \to B^{n_{\ell}}$ is a $\frac{1}{1 - s}$-Lipschitz homeomorphism \cite[p.~201]{SW}. 
(b) If ${\boldsymbol t}_{\ell} \in \sigma$, then $\lambda_{\ell}({\boldsymbol t}_{\ell})$ is a PL-homeomorphism. 
(c) $\lambda_{\ell}({\boldsymbol 0}) = id_{B^{n_{\ell}}}$.
\end{proof}

If we replace $s^\infty$ by $s$ and omit the indices $k$ and $\ell$ in Lemmas 4.4 and 4.5, 
then we obtain the corresponding statements for the case (I). 

In the proofs of Propositions 4.1~(i) and 4.2~(i) we apply the construction (I). 

\begin{proof}[Proof of Proposition 4.1 (i)] 
We choose any homeomorphism $\chi : (s^2, \Sigma^2, \sigma^2) \cong (s, \Sigma, \sigma)$ 
with $\chi({\boldsymbol 0}, {\boldsymbol 0}) = {\boldsymbol 0}$ and 
a similar affine embedding $\phi : B^n \to {\rm Int}\,C$, where 
$C$ is a principal simplex of in $X \setminus A$ of $n \equiv {\rm dim}\,C \geq 1$. 
For $Y = X$, $\phi$ and $\chi$ we obtain the stability homeomorphism $\Phi : {\cal F} \times s \to {\cal F}$. 

The following pairs $({\cal F}_1, S_1)$ (${\cal G}_1 = {\cal F}_1$) satisfy the conditions (a) -- (d) in Lemma 4.4:
\[ ({\mathcal H}_A(X), s), 
({\mathcal H}^{\rm loc\,LIP}_A(X), \Sigma), 
({\mathcal H}^{\rm LIP}_A(X), \Sigma), 
({\mathcal H}^{\rm LIP, c}_A(X), \Sigma), 
({\mathcal H}^{\rm PL}_A(X), \sigma), 
({\mathcal H}^{\rm PL, c}_A(X), \sigma) \]
\noindent In fact, the conditions (a) and (b) follow from Lemma 4.5~(i) and (ii) respectively and (c) is obvious.  
Thus $\Phi$ induces the required stability homeomorphism of tuples. 
Since $\Phi(id_X, {\boldsymbol 0}) = id_X$, it follows that $\Phi(({\cal F}_1)_0 \times S_1) = ({\cal F}_1)_0$ for each pair $({\cal F}_1, S_1)$.
\end{proof}

\begin{proof}[Proof of Proposition 4.2 (i)]
Choose a PL-homeomorphism $h : M \to Y$ onto a Euclidean polyhedron and let $A = h(X)$. 
By Lemma 4.2 $h$ induces a homeomorphism of tuples, defined by $h^\#(f) = h f h^{-1}$:
\begin{equation*}
\begin{split}
h^\# : \ & ({\cal H}_X(M), {\cal H}^{\rm loc\,LIP}_X(M), {\cal H}^{\rm LIP}_X(M), {\cal H}^{\rm LIP, c}_X(M), 
{\mathcal H}^{\rm PL}_X(M), {\cal H}^{\rm PL, c}_X(M)) \\
& \hspace{38mm} \cong \ ({\cal H}_A(Y), {\cal H}^{\rm loc\,LIP}_A(Y), {\cal F}_1, {\cal H}^{\rm LIP, c}_A(Y), {\mathcal H}^{\rm PL}_A(Y), {\cal H}^{\rm PL, c}_A(Y))
\end{split}
\end{equation*}
where ${\cal F}_1 = h^\#({\cal H}^{\rm LIP}_X(M))$. 
It suffices to show that the latter tuple is $(s, \Sigma, \Sigma, \Sigma, \sigma, \sigma)$-stable. 

Let $F$, $G$ and $\Phi$ be the maps given in the proof of Proposition 4.1 (i) for $(Y, A)$.   
Then $F({\cal H}^{\rm loc\,LIP}_A(Y)) \subset \Sigma$ and $G(\Sigma) \subset {\cal H}^{\rm LIP, c}_A(Y)$. 
Since ${\cal H}^{\rm loc\,LIP}_A(Y) \supset {\cal F}_1 \supset {\cal H}^{\rm LIP, c}_A(X)$, 
the data ${\cal F}_1$, $\Sigma$ and ${\cal G}_1 = {\cal H}^{\rm LIP, c}_A(Y)$ satisfy the conditions (a) -- (d) in Lemma 4.4, 
and this implies that $\Phi({\cal F}_1 \times \Sigma) = {\cal F}_1$. This completes the proof. 
\end{proof}

In the proofs of Propositions 4.1 (ii) and 4.2~(iii) we apply the construction (II) using the homeomorphism
\begin{equation*}
\begin{split}
\chi &: ((s^{\infty})^2, (\Sigma^{\infty})^2, (s^\infty_b)^2, (\Sigma^{\infty}_f)^2, (\sigma^{\infty})^2, (\sigma^{\infty}_f)^2) 
\cong (s^{\infty}, \Sigma^{\infty}, s^\infty_b, \Sigma^{\infty}_f, \sigma^{\infty}, \sigma^{\infty}_f) \\
& \hspace{25mm} \chi(({\boldsymbol t}_1, {\boldsymbol t}_2, \cdots), ({\boldsymbol s}_1, {\boldsymbol s}_2, \cdots)) 
= ({\boldsymbol t}_1, {\boldsymbol s}_1, {\boldsymbol t}_2, {\boldsymbol s}_2, \cdots)
\end{split}
\end{equation*}

\begin{proof}[Proof of Proposition 4.1 (ii)]
We choose (a) a triangulation $T$ of $X$ with ${\rm diam}\,C \leq 1$ for each simplex $C$ of $T$, 
(b) a sequence of principal simplices $\{ C_k \}_{k=1}^\infty$ of $T$ 
such that $n_k \equiv \dim \, C_k \geq 1$, $C_k \subset X \setminus A$ and $d(C_k, C_\ell) \geq 1$ ($k \neq \ell$), and 
(c) similar affine embeddings $\phi_k : B^{n_k} \to {\rm Int}\,C_k$ such that
${\rm diam}\,\phi_k(B^{n_k}) < \varepsilon_k$ and $N_X(\phi_k(B^{n_k}), \varepsilon_k) \subset {\rm Int}\,C_k$ for some $\varepsilon_k \in (0, 1)$.

For $Y = X$, $\phi_k$ and $\chi$ we obtain the stability homeomorphism $\Phi : {\cal F} \times s^\infty \to {\cal F}$. 

By the choice of $\phi_k$ and Lemma 4.5 the following pairs $({\cal F}_1, S_1)$ (${\cal G}_1 = {\cal F}_1$) satisfy the conditions (a) -- (d) in Lemma 4.4:
\[ ({\mathcal H}_A(X), s^{\infty}), 
({\mathcal H}^{\rm loc\,LIP}_A(X), \Sigma^{\infty}), 
({\mathcal H}^{\rm LIP}_A(X), s^\infty_b), 
({\mathcal H}^{\rm LIP, c}_A(X), \Sigma^{\infty}_f), 
({\mathcal H}^{\rm PL}_A(X), \sigma^{\infty}), 
({\mathcal H}^{\rm PL, c}_A(X), \sigma^{\infty}_f) \]
\noindent 
Note that for ${\boldsymbol t} = ({\boldsymbol t}_k)_k \in s^\infty$, 
(a) if ${\boldsymbol t}_{\ell} \in [-s, s]^\infty$ ($0 \leq s < 1$), 
then $G({\boldsymbol t}) : \phi_\ell(B^{n_\ell}) \to \phi_\ell(B^{n_\ell})$ and also 
$G({\boldsymbol t}) : C_{\ell} \to C_{\ell}$ are $1/(1-s)$-Lipschitz homeomorphisms cf. \cite[Lemma 1.4]{SW}, and 
(b) if ${\boldsymbol t}_k \in [-s, s]^\infty$ for each $k$, then $G({\boldsymbol t}) : X \to X$ is a $K$-Lipschitz homeomorphism for $K = \max\{1/(1-s), 3\}$ by (a) and the choice of $C_k$ and $\varepsilon_k$. 

Thus $\Phi$ induces the required stability homeomorphism of tuples. 
Since $\Phi(id_X, {\boldsymbol 0}) = id_X$, the homeomorphism $\Phi$ preserves the identity components of the homeomorphism groups.
\end{proof}

Proposition 4.2 (ii) follows from Proposition 4.1 (ii) and Lemma 4.2.

\begin{proof}[Proof of Proposition 4.2 (iii)]
We may assume that $p_k \in {\rm Int}\,M \setminus X$. 
For each $k \geq 1$ we choose a $\varepsilon_k$-open neighborhood $U_k$ of the origin $0$ in $T_{p_k}M$
such that the exponential map $\exp_{p_k}$ maps $U_k$ diffeomorphically onto a small open neighborhood $V_k$ of $p_k$ in ${\rm Int}\,M \setminus X$   
with ${\rm diam}\,V_k \leq \varepsilon/3$. 
Since $\exp_{p_k}$ is isometric at the origin of $T_{p_k}M$, if $\varepsilon_k$ is so small, then 
$\exp_{p_k} : U_k \to V_k$ is a 2-Lipschitz homeomorphism. 
We choose any $c_k$-similar linear isomorphism $\theta_k : {\Bbb R}^n \to T_{p_k}M$ ($c_k > 0$ small) and 
consider the embeddings $\phi_k = \exp_{p_k} \theta_k : B^n \to M$. 
If $c_k$ is small, then ${\rm diam}\,\phi_k(B^n) < \delta_k$ and $N_M(\phi_k(B^n), \delta_k) \subset V_k$ for some $\delta_k \in (0, \varepsilon/3)$.
Note that $\phi_k : [-1, 1] \to M$, $\phi_k(t) = \exp_{p_k}\,tv$ ($v = \theta_k(1)$, $\|v\|_{p_k} = c_k$), 
is a $c_k$-similar homeomorphism onto a shortest geodesic $\phi_k([-1, 1])$. 
For $Y = X$, $\phi_k$ and $\chi$ we obtain the stability homeomorphism $\Phi : {\cal F} \times s^\infty \to {\cal F}$. 

By the choice of $\phi_k$ and Lemma 4.5 the following pairs $({\cal F}_1, S_1)$ (${\cal G}_1 = {\cal F}_1$) satisfy the conditions (a) -- (d) in Lemma 4.4:
\[ 
({\mathcal H}_X(M), s^{\infty}), 
({\mathcal H}^{\rm loc\,LIP}_X(M), \Sigma^{\infty}), 
({\mathcal H}^{\rm LIP}_X(M), s^\infty_b), 
({\mathcal H}^{\rm LIP, c}_X(M), \Sigma^{\infty}_f)
\]
\noindent Note that for ${\boldsymbol t} = ({\boldsymbol t}_k)_k \in s^\infty$, 
(a) if ${\boldsymbol t}_{\ell} \in [-s, s]^\infty$ ($0 \leq s < 1$), 
then $G({\boldsymbol t}) : V_{\ell}\to V_{\ell}$ is a $4/(1-s)$-Lipschitz homeomorphism, and 
(b) if ${\boldsymbol t}_k \in [-s, s]^\infty$ for each $k$, 
then $G({\boldsymbol t}) : X \to X$ is a $K$-Lipschitz homeomorphism for $K = \max\{4/(1-s), 3\}$. 
For the claim (a), we observe that 
$G({\boldsymbol t})|_{V_\ell} = \exp_{p_\ell} \, \theta_\ell \, \widetilde{\lambda_\ell({\boldsymbol t}_\ell)} \, (\theta_\ell)^{-1} \, (\exp_{p_\ell})^{-1}$, 
where $\exp_{p_\ell} : U_\ell \to V_\ell$ is a 2-Lipschitz homeomorphism, 
$\theta_\ell : W_\ell \equiv \theta_\ell^{-1}(U_\ell) \to U_\ell$ is $c_\ell$-similar, 
and $\widetilde{\lambda_\ell({\boldsymbol t}_\ell)} : W_\ell \to W_\ell$ is the extension of $\lambda_\ell({\boldsymbol t}_\ell) : B^n \to B^n$  
by the identity, which is a $\frac{1}{1-s}$-Lipschitz homeomorphism since $\lambda_\ell({\boldsymbol t}_\ell)$ is so.
The claim (b) follows from (a) and the choice of $\varepsilon$, $V_\ell$ and $\delta_\ell$.

Therefore $\Phi$ induces the required stability homeomorphism of tuples. 
Since $\Phi(id_M, {\boldsymbol 0}) = id_M$, the homeomorphism $\Phi$ preserves the identity components of the homeomorphism groups.
\end{proof}

\section{$(s^{\infty}, \Sigma^{\infty}, \Sigma_f^{\infty}, \sigma_f^{\infty})$-manifolds} 

\subsection{Characterization of $(s, S_1, \cdots, S_l)$-manifolds} \mbox{}
\vskip 2mm
First we recall a general characterization of $(s, S_1, \cdots, S_l)$-manifold based upon the stability property \cite[Section 2.2]{Ya1}. 
We represent $s$ as $s = \prod_{k \in {\Bbb N}} {\Bbb R}$ (${\Bbb N} = \{ 1, 2, \cdots \}$).  
For $A \subset {\Bbb N}$ we set $c(A) = {\Bbb N} \setminus A$, $s(A) = \prod_{k \in A} \, (-\infty, \infty)$ 
and let $\pi_A : s \to s(A)$ denote the projection. For a subset $S$ of $s$ let $S(A) = \pi_A(S) \subset s(A)$. 

We assume that the model tuple $(s, S_1, \cdots, S_l)$ satisfies the following condition (\#): 
\vskip 2mm
\noindent (\#1) each $S_i$ is a linear subspace of $s$ and $S_1$ is a $\sigma$ $Z$-set of $S_1$ itself, \\
(\#2) $S_l$ has the h.n. complement in $s$, \\ 
(\#3) there exists a sequence $A_n$ ($n \geq 1$) of disjoint infinite subsets of ${\Bbb N}$ such that 
for each $i = 1, \cdots, l$ and $n \geq 1$ 
(a) ${\rm min} \, A_n > n$, 
(b) $S_i = S_i(A_n) \times S_i(c(A_n))$ and 
(c) $(s(A_n), S_1(A_n), \cdots, S_l(A_n)) \cong (s, S_1, \cdots, S_l)$.  
\vskip 2mm
Let ${\mathcal M}(s, S_1, \cdots, S_l)$ denote the class of $(l+1)$-tuples $(X, X_1, \cdots, X_l)$ 
which has a closed embedding $h : X \to s$ such that $h^{-1}(S_i) = X_i \ (1 \leq i \leq l)$. 

\begin{theorem}
A tuple $(X, X_1, \cdots, X_l)$ is a $(s, S_1, \cdots, S_l)$-manifold iff \\ 
(i) $X$ is a separable completely metrizable ANR, \\
(ii) $X_l$ has the h.n. complement in $X$, \\ 
(iii) $(X, X_1, \cdots, X_l) \in {\mathcal M}(s, S_1, \cdots, S_l)$, \\ 
(iv)(Stability) $(X, X_1, \cdots, X_l)$ is $(s, S_1, \cdots, S_l)$-stable.
\end{theorem}

\begin{proposition}(Homotopy Invariance) 
Suppose $(X, X_1, \cdots, X_l)$ and $(Y, Y_1, \cdots, Y_l)$ are $(s, S_1, \cdots, S_l)$-manifolds. 
Then $(X, X_1, \cdots, X_l) \cong (Y, Y_1, \cdots, Y_l)$ iff $X \simeq Y$.
\end{proposition} 

\subsection{Characterization of $(s^{\infty}, \Sigma^{\infty}, \Sigma_f^{\infty}, \sigma_f^{\infty})$-manifolds} \mbox{} 
\vskip 2mm 
This subsection is devoted to the verification of Theorem 2.1 and Proposition 2.1. 
These statements follow from Theorem 5.1, Proposition 5.1 and the next lemma. 

\begin{lemma} 
{\rm (1)} The quadruple $(s^{\infty}, \Sigma^{\infty}, \Sigma_f^{\infty}, \sigma_f^{\infty})$ satisfies the condition (\#). 

{\rm (2)} The class ${\mathcal M}(s^{\infty}, \Sigma^{\infty}, \Sigma_f^{\infty}, \sigma_f^{\infty})$ coincides with 
the class ${\mathcal M}$ of quadruples $(X, X_1, X_2, X_3)$ such that 
$X$ is completely metrizable, $X_1$ is $F_{\sigma \delta}$ in $X$, $X_2$ is $\sigma$-compact and $X_3$ is $\sigma$-fd-compact. 
\end{lemma}

\begin{proof}
(1) (\#1) $\Sigma^{\infty}$ is a $\sigma$ $Z$-set of $\Sigma^{\infty}$ itself 
since $[-n, n]^\infty$ ($n \geq 1$) are $Z$-sets of $\Sigma$ and $\Sigma^{\infty} = \Sigma \times \Sigma^{\infty} = \cup_n ([-n, n]^\infty \times \Sigma^{\infty})$. 

(\#2) Note that $(s^{\infty}, \sigma_f^{\infty}) \cong (s, \sigma)$. 

(\#3) We identify $s$ with 
$s^{\infty} = \prod_{k \in {\Bbb N}} s = \prod_{(k, i) \in {\Bbb N} \times {\Bbb N}} {\Bbb R}$ 
by taking any linear ordering of ${\Bbb N} \times {\Bbb N}$. 
Take an appropriate sequence $B_n \ (n \geq 1)$ of disjoint infinite subsets of ${\Bbb N}$ so that $A_n \equiv {\Bbb N} \times B_n$ ($n \geq 1)$ satisfies the condition (a). It follows that (c) $(s^{\infty}(A_n), \Sigma^{\infty}(A_n), \Sigma_f^{\infty}(A_n), \sigma_f^{\infty}(A_n)) = (s(B_n)^{\infty}, \Sigma(B_n)^{\infty}, \Sigma(B_n)_f^{\infty}, \sigma(B_n)_f^{\infty}) \cong (s^{\infty}, \Sigma^{\infty}, \Sigma_f^{\infty}, \sigma_f^{\infty})$ 
and (b) $(\Sigma^{\infty}, \Sigma_f^{\infty}, \sigma_f^{\infty}) = (\Sigma(B_n)^{\infty} \times \Sigma(c(B_n))^{\infty}, \Sigma(B_n)_f^{\infty} \times \Sigma(c(B_n))_f^{\infty}, \sigma(B_n)_f^{\infty} \times \sigma(c(B_n))_f^{\infty})$. 

(2) The assertion follows from the next Lemma 5.2. The required closed embedding $f : X \to s^{\infty}$ is defined by $f = \chi (H, \phi)$, where $\chi : ((s^{\infty})^2$, $(\Sigma^{\infty})^2$, $(\Sigma_f^{\infty})^2, (\sigma_f^{\infty})^2) \cong (s^{\infty}, \Sigma^{\infty}, \Sigma_f^{\infty}, \sigma_f^{\infty})$ is a homeomorphism of quadruples. 
\end{proof}

\begin{lemma}
{\rm (1)} If $(X, X_1, X_2, X_3)$ is a quadruple such that $X_1$ is $F_{\sigma \delta}$ and $X_2$ and $X_3$ are $F_{\sigma}$ in $X$, then there exists a map $\phi : X \to s^{\infty}$ such that $\phi^{-1}(\Sigma^{\infty}) = X_1$, $\phi^{-1}(\Sigma_f^{\infty}) = X_2$ and $\phi^{-1}(\sigma_f^{\infty}) = X_3$.  \\ 
{\rm (2)} For each $(X, X_1, X_2, X_3) \in {\mathcal M}$ there exists a closed embedding $H : X \to s^{\infty}$ such that $H(X) \subset \Sigma^{\infty}$, $H(X_2) \subset \Sigma^{\infty}_f$ and $H(X_3) \subset \sigma^{\infty}_f$. 
\end{lemma}

\begin{proof} 
The proof is a modification of the argument of \cite[\S3.1]{Ya1}. 

(1) We can write $X_1 = \cap_n  A_n$, where $A_n$ is a $F_{\sigma}$ subset of $X$, and $X_2 = \cup_n \, B_n$, where $B_n$ is a closed subset of $X$ and $B_n \subset B_{n+1}$. By \cite[Lemma 3.3 (i)]{Ya1} there exist 
(a) a map $\alpha : X \to (-1, 1)^\infty \subset s$ such that $\alpha^{-1}((-1, 1)^\infty_f) = X_3$ and 
(b) maps $\beta_n : X \to s$ such that $\beta_n^{-1}(\Sigma) = A_n$, $\beta_n(A_n) \subset \sigma$ and $\beta_n^{-1}(0) = B_n$. 
The required map $\phi : X \to s^{\infty}$ is defined by $\phi = (\alpha, \beta_1, \beta_2, \cdots)$.

(2) Since $X_2$ is $\sigma$-compact, there exists a closed embedding $f : X \to s$ such that $f(X_2) \subset \Sigma$ \cite{Cp}. Since $\Sigma = \cup_n [-n, n]^{\infty}$, we can write $X_2 = \cup_{n=1}^{\infty} C_n$ so that each $C_n$ is compact, $C_n \subset C_{n+1}$ and $f(C_n) \subset [-n, n]^{\infty}$. Since (a) $[-(n-1), n-1]^{\infty}$ is a $Z$-set of $[-n, n]^{\infty}$, (b) $[-n, n]^{\infty}$ is a $Z$-set of $s$ and (c) $\sigma$ has the h.n. complement in $s$, 
there exists a map $f_n : X \to \Sigma$ such that 
$f_n(X \setminus C_n) \subset \sigma \setminus [-n, n]^{\infty}$, 
$f_n|_{C_n} : C_n \to [-n, n]^{\infty}$ is an embedding, 
$f_n(C_{n-1}) \subset [-(n-1), n-1]^{\infty}$, $f_n(C_n \setminus C_{n-1}) \subset [-n, n]^{\infty} \setminus [-(n-1), n-1]^{\infty}$, 
$f_n(X_3 \cap C_n) \subset [-n, n]^\infty_f$ and 
$d(f_n, f) < \frac{1}{n}$, where $C_0 = \emptyset$ and $d$ is a fixed complete metric on $s$. 
By \cite[Lemma 3.4 (i)]{Ya1} the map $F = (f_n)_n : X \to s^{\infty}$ is a closed embedding such that $F(X) \subset \Sigma^{\infty}$, $\displaystyle F(C_n) \subset \Sigma^{n-1} \times [-n, n]^\infty \times \prod_{m \geq n+1} [-(m-1), m-1]^\infty$ and $F(X_3) \subset \sigma^{\infty}$. 

We regard as $Q = [-\infty, \infty]^\infty$. For each $n \geq 1$, there exists a map $p_n : Q \to Q$ such that (a) $p_n^{-1}(0) = [-(n-1), n-1]^{\infty}$, (b) $p_n$ maps $Q \setminus [-(n-1), n-1]^{\infty}$ homeomorphically onto $Q \setminus 0$, (c) $p_n^{-1}(s) = s$ and (d) if $p_n(x) = y$ then $|y_i| \leq |x_i| \ (i \geq 1)$ so that $p_n(\Sigma) \subset \Sigma$ and $p_n(\sigma) \subset \sigma$. This follows from an easy shrinking argument based on the next claim. 

\begin{claim}
For any neighborhoods $U$ of $[-n, n]^{\infty}$ and $V$ of $0$ in $Q$, there exists a homeomorphism $h : Q \to Q$ such that $h = id$ on $Q \setminus U$, $h(Q \setminus s) = Q \setminus s$, $h([-n, n]^{\infty}) \subset V$ and if $h(x) = y$, then $|y_i| \leq |x_i|$ for each $i \geq 1$. 
\end{claim}

The homeomorphism $h$ is obtained as $h = h_m \times id$, where $m$ is sufficiently large and $h_m$ shrinks $[-n, n]^m$ radially towards $0$ in $[-\infty, \infty]^m$. By the repeated application of Claim, we can construct a sequence of homeomorphisms of $Q$ which converges to the desired map $p_n : Q \to Q$. 

The map $G = \prod_n p_n : s^{\infty} \to s^{\infty}$ is a closed map and $G(\Sigma^{\infty}) \subset \Sigma^{\infty}$, $G(\sigma^{\infty}) \subset \sigma^{\infty}$. Hence the map $H = GF = (p_nf_n)_n : X \to s^{\infty}$ is a closed map and $H(X) \subset \Sigma^{\infty}$. It remains to show that (i) $H(X_2) \subset \Sigma^{\infty}_f$, (ii) $H(X_3) \subset \sigma^{\infty}_f$ and (iii) $H$ is injective. (i) If $x \in X_2$, then $x \in C_n$ for some $n$ and it follows that $f_m(x) \in [-(m-1), m-1]^{\infty}$ and hence $p_mf_m(x) = 0$ for each $m \geq n + 1$. 
(ii) $H(X_3) \subset \Sigma^\infty_f \cap \sigma^\infty = \sigma^\infty_f$. 
(iii) Suppose $x, y \in X$ and $x \neq y$. 
If $x$, $y \in X_2$ then we may assume that $x \in C_n \setminus C_{n-1}$ and $y \in C_n$. 
The choice of $f_n$ implies that $f_n(x) \neq f_n(y)$ and $f_n(x) \not\in [-(n-1), n-1]^{\infty}$.
By the choice $p_n$ we have $p_nf_n(x) \neq p_nf_n(y)$ and so $H(x) \neq H(y)$.
If $x \not\in X_2$, then since $F$ is injective, we have $f_n(x) \neq f_n(y)$ for some $n \geq 1$, 
and since $x \not\in C_n$ we have $f_n(x) \not\in [-n, n]^{\infty}$. 
Again, by the choice $p_n$ we have $p_nf_n(x) \neq p_nf_n(y)$ and $H(x) \neq H(y)$.
This completes the proof.    
\end{proof}        

\subsection{Characterization of $(s^{\infty}, \Sigma^{\infty}, s^{\infty}_b, \sigma_f^{\infty})$-manifolds} \mbox{} 
\vskip 2mm 
The next proposition clarifies the relation between $s^{\infty}_b$ and $\Sigma_f^{\infty}$ in $s^{\infty}$. 

\begin{proposition} 
$(s^{\infty}, \Sigma^{\infty}, \Sigma_f^{\infty}, \sigma_f^{\infty}) \cong (s^{\infty}, \Sigma^{\infty}, s^{\infty}_b, \sigma_f^{\infty})$.  
\end{proposition}

To prove this statement we need a similar characterization of $(s^{\infty}, \Sigma^{\infty}, s^{\infty}_b, \sigma_f^{\infty})$-manifolds. 

\begin{proposition} 
A quadruple $(X, X_1, X_2, X_3)$ is a $(s^{\infty}, \Sigma^{\infty}, s^{\infty}_b, \sigma_f^{\infty})$-manifold iff \\
(i) $X$ is a completely metrizable ANR, \\
(ii) $X_3$ has the h.n. complement in $X$, \\
(iii) $X_1$ is $F_{\sigma \delta}$ in $X$, $X_2$ is $\sigma$-compact, $X_3$ is $\sigma$-fd-compact, \\
(iv) $(X, X_1, X_2, X_3)$ is $(s^{\infty}, \Sigma^{\infty}, s^{\infty}_b, \sigma_f^{\infty})$-stable. 
\end{proposition} 

\begin{proposition}
Suppose $(X, X_1, X_2, X_3)$ and $(Y, Y_1, Y_2, Y_3)$ are $(s^{\infty}, \Sigma^{\infty}, s^{\infty}_b, \sigma_f^{\infty})$-manifolds. 
Then $(X, X_1, X_2, X_3) \cong (Y, Y_1, Y_2, Y_3)$ iff $X \simeq Y$.
\end{proposition} 

The quadruple $(s^{\infty}, \Sigma^{\infty}, s_b^{\infty}, \sigma_f^{\infty})$ satisfies the condition (\#). 
This is verified by the same observation as in Lemma 5.1. Therefore Propositions 5.3 and 5.4 follow from Theorem 5.1, Proposition 5.1 and the following lemma. 

\begin{lemma} 
The class ${\mathcal M}(s^{\infty}, \Sigma^{\infty}, s_b^{\infty}, \sigma_f^{\infty})$ coincides with the class ${\mathcal M}$ of quadruples $(X, X_1, X_2, X_3)$ such that $X$ is completely metrizable, $X_1$ is $F_{\sigma \delta}$ in $X$, $X_2$ is $\sigma$-compact and $X_3$ is $\sigma$-fd-compact. 
\end{lemma}

\begin{proof}
Let $(X, X_1, X_2, X_3) \in {\mathcal M}$. By Lemma 5.1(2) there exists a closed embedding $h : X \to s^\infty$ such that 
$h^{-1}(\Sigma^\infty) = X_1$, $h^{-1}(\Sigma^\infty_f) = X_2$ and $h^{-1}(\sigma^\infty_f) = X_3$. 
Note that $h(X_2) \subset \Sigma^\infty_f \subset s^\infty_b$. There exists a map $\alpha : X \to s$ 
such that $\alpha^{-1}(\Sigma) = X_2$ and $\alpha^{-1}(\sigma) = X_3$ \cite[Lemma 3.3\,(ii)]{Ya1}. 
The required closed embedding $f : X \to s^{\infty}$ is defined by $f = \psi (\alpha, h)$, 
where the homeomorphism $\psi : (s \times s^{\infty}, s \times \Sigma^{\infty}, \Sigma \times s_b^{\infty}) \cong (s^{\infty}, \Sigma^{\infty}, s_b^{\infty})$ is defined by $\psi((x_i)_i, ({\boldsymbol y}_k)_k) = ((x_k, {\boldsymbol y}_k))_k$. 
\end{proof}

\vspace{2mm} 

\begin{proof}[Proof of Proposition 5.2.] 
Since $((s^{\infty})^2, (\Sigma^{\infty})^2, (s_b^{\infty})^2, (\Sigma^{\infty}_f)^2, (\sigma_f^{\infty})^2) \cong (s^{\infty}, \Sigma^{\infty}, s_b^{\infty}, \Sigma_f^{\infty}, \sigma_f^{\infty})$, 
the product $((s^{\infty})^2, (\Sigma^{\infty})^2, \Sigma_f^{\infty} \times s_b^{\infty}, (\sigma_f^{\infty})^2)$ satisfies both $(s^{\infty}, \Sigma^{\infty}, \Sigma_f^{\infty}, \sigma_f^{\infty})$- and $(s^{\infty}, \Sigma^{\infty}, s_b^{\infty}, \sigma_f^{\infty})$-stability. Hence by Theorem 2.1 and Proposition 5.3 $((s^{\infty})^2, (\Sigma^{\infty})^2, \Sigma_f^{\infty} \times s_b^{\infty}, (\sigma_f^{\infty})^2)$ is both an $(s^{\infty}, \Sigma^{\infty}, \Sigma_f^{\infty}, \sigma_f^{\infty})$- and an $(s^{\infty}, \Sigma^{\infty}, s_b^{\infty}, \sigma_f^{\infty})$-manifold. Thus by Propositions 2.1 and 5.4 we have 
\[ (s^{\infty}, \Sigma^{\infty}, \Sigma_f^{\infty}, \sigma_f^{\infty}) \cong ((s^{\infty})^2, (\Sigma^{\infty})^2, \Sigma_f^{\infty} \times s_b^{\infty}, (\sigma_f^{\infty})^2) \cong (s^{\infty}, \Sigma^{\infty}, s_b^{\infty}, \sigma_f^{\infty}).\]  
\end{proof} 

\section{Some problems}

We conclude this paper with some problems. 

\begin{problem}
Suppose $M$ is a connected noncompact complete Riemannian 2-manifold with a $C^1$-triangulation and $X$ is a compact subpolyhedron of $M$.
Is the quadruple 
\[ ({\cal H}_X(M)_0, {\cal H}_X^{\rm loc \, LIP}(M)_0, {\cal H}_X^{\rm LIP}(M)_0, {\cal H}_X^{\rm PL, c}(M)_0) \]
a $(s^{\infty},\Sigma^{\infty},\Sigma^{\infty}_f, \sigma^{\infty}_f)$-manifold ?
\end{problem}

As noted in Section 4.1, if every complete Riemannian manifold admits a Lipschitz PL-embedding into some Euclidean space, 
then this follows from the corresponding statement for Euclidean PL 2-manifolds.

\begin{problem} 
Suppose $M$ is a Riemann surface with a QC-triangulation and $X$ is a compact subpolyhedron of $M$. \\
(i) When $M$ is compact, is the triple $({\cal H}_X(M), {\cal H}_X^{\rm QC}(M), {\cal H}_X^{\rm PL}(M))$ a $(s, \Sigma, \sigma)$-manifold ? \\
(ii) When M is connected and noncompact, is the quadruple 
$({\cal H}_X(M)_0, {\cal H}_X^{\rm loc \, QC}(M)_0, {\cal H}_X^{\rm QC}(M)_0, {\cal H}_X^{\rm PL, c}(M)_0)$ 
a $(s^{\infty},\Sigma^{\infty},\Sigma^{\infty}_f, \sigma^{\infty}_f)$-manifold ?
\end{problem}

In \cite{Ya4} we have proved the corresponding statement for the pair $({\cal H}(M)_0, {\cal H}^{\rm QC}(M)_0)$. 
However, our stability argument does not enable us to include the subgroup ${\cal H}_X^{\rm PL, c}(M)_0$

Lipschitz homeomorphisms are H\"older continuous and QC-homeomorphisms are locally H\"older continuous. 
We expect that this observation leads to a systematic treatment of these groups.

\begin{problem} 
Is there a unified approach to treat the groups of H\"older, QC, Lipschitz and PL homeomorphisms ?
\end{problem}

\end{document}